\documentclass{amsart}

\usepackage{amsmath,amsthm,amsfonts,amscd,amssymb}

\newtheorem{thm}{Theorem}[section]
\newtheorem{cor}[thm]{Corollary}
\newtheorem{lem}[thm]{Lemma}

\theoremstyle{definition}

\theoremstyle{remark}

\begin{document}

\title{On effective Witt decomposition and Cartan-Dieudonn{\'e} theorem}
\author{Lenny Fukshansky}

\address{Department of Mathematics, Mailstop 3368, Texas A\&M University, College Station, Texas 77843-3368}
\email{lenny@math.tamu.edu}

\begin{abstract}
Let $K$ be a number field, and let $F$ be a symmetric bilinear form in $2N$ variables over $K$. Let $Z$ be a subspace of $K^N$. A classical theorem of Witt states that the bilinear space $(Z,F)$ can be decomposed into an orthogonal sum of hyperbolic planes, singular, and anisotropic components. We prove the existence of such a decomposition of small height, where all bounds on height are explicit in terms of heights of $F$ and $Z$. We also prove a special version of Siegel's Lemma for a bilinear space, which provides a small-height orthogonal decomposition into one-dimensional subspaces. Finally, we prove an effective version of Cartan-Dieudonn{\'e} theorem. Namely, we show that every isometry $\sigma$ of a regular bilinear space $(Z,F)$ can be represented as a product of reflections of bounded heights with an explicit bound on heights in terms of heights of $F$, $Z$, and $\sigma$.
\end{abstract}

\maketitle

\def\A{{\mathcal A}}
\def\B{{\mathcal B}}
\def\C{{\mathcal C}}
\def\D{{\mathcal D}}
\def\F{{\mathcal F}}
\def\x{{\mathcal H}}
\def\I{{\mathcal I}}
\def\J{{\mathcal J}}
\def\K{{\mathcal K}}
\def\L{{\mathcal L}}
\def\M{{\mathcal M}}
\def\R{{\mathcal R}}
\def\s{{\mathcal S}}
\def\V{{\mathcal V}}
\def\X{{\mathcal X}}
\def\Y{{\mathcal Y}}
\def\H{{\mathcal H}}
\def\OO{{\mathcal O}}
\def\cee{{\mathbb C}}
\def\pee{{\mathbb P}}
\def\que{{\mathbb Q}}
\def\real{{\mathbb R}}
\def\zed{{\mathbb Z}}
\def\hyp{{\mathbb H}}
\def\qbar{{\overline{\mathbb Q}}}
\def\eps{{\varepsilon}}
\def\ahat{{\hat \alpha}}
\def\bhat{{\hat \beta}}
\def\hs{{\hat \sigma}}
\def\gt{{\tilde \gamma}}
\def\h{{\tfrac12}}
\def\be{{\boldsymbol e}}
\def\bei{{\boldsymbol e_i}}
\def\baf{{\boldsymbol f}}
\def\baa{{\boldsymbol a}}
\def\bc{{\boldsymbol c}}
\def\bm{{\boldsymbol m}}
\def\bk{{\boldsymbol k}}
\def\bi{{\boldsymbol i}}
\def\bl{{\boldsymbol l}}
\def\bq{{\boldsymbol q}}
\def\bu{{\boldsymbol u}}
\def\bt{{\boldsymbol t}}
\def\bs{{\boldsymbol s}}
\def\bv{{\boldsymbol v}}
\def\bw{{\boldsymbol w}}
\def\bx{{\boldsymbol x}}
\def\bX{{\boldsymbol X}}
\def\bz{{\boldsymbol z}}
\def\bwy{{\boldsymbol y}}
\def\bY{{\boldsymbol Y}}
\def\bL{{\boldsymbol L}}
\def\ba{{\boldsymbol\alpha}}
\def\bb{{\boldsymbol\beta}}
\def\bet{{\boldsymbol\eta}}
\def\bxi{{\boldsymbol\xi}}
\def\bo{{\boldkey 0}}
\def\bol{{\boldkey 1}_L}
\def\ep{\varepsilon}
\def\p{\boldsymbol\varphi}
\def\q{\boldsymbol\psi}
\def\rank{\operatorname{rank}}
\def\aut{\operatorname{Aut}}
\def\lcm{\operatorname{lcm}}
\def\sgn{\operatorname{sgn}}
\def\spn{\operatorname{span}}
\def\md{\operatorname{mod}}
\def\Norm{\operatorname{Norm}}
\def\dim{\operatorname{dim}}
\def\det{\operatorname{det}}
\def\Vol{\operatorname{Vol}}
\def\rk{\operatorname{rk}}

\section{Introduction and notation}

Let $K$ be a number field, $N>1$ an integer. Let
$$F(\bX,\bY) = \sum_{i=1}^N \sum_{j=1}^N f_{ij} X_i Y_j,$$
be a symmetric bilinear form in $2N$ variables with coefficients $f_{ij} = f_{ji}$ in $K$. We will write $F(\bX) = F(\bX,\bX)$ for the associated quadratic form in $N$ variables, and will also use $F$ to denote the symmetric $N \times N$ matrix $(f_{ij})_{1 \leq i,j \leq N}$. Let $Z \subseteq K^N$ be an $L$-dimensional subspace, $2 \leq L \leq N$, then $F$ is also defined on $Z$, and we write $(Z,F)$ for the bilinear space. Let $M$ be the Witt index of $(Z,F)$. With this basic notation we can recall the classical Witt decomposition theorem. We give a brief overview of required definitions and basic results on bilinear spaces in section~3. 

\begin{thm} \label{witt} Suppose that $(Z,F)$ is a bilinear space as above. Then there exists an orthogonal decomposition of $(Z,F)$ of the form
\begin{equation}
\label{decompose}
Z = Z^{\perp} \perp \hyp_1 \perp\ ...\ \perp \hyp_M \perp V,
\end{equation}
where $Z^{\perp} = \{ \bx \in Z : F(\bx,\bz) = 0\ \forall\ \bz \in Z \}$ is the singular component, $\hyp_i$ are hyperbolic planes, and $V$ is anisotropic component, which is uniquely determined up to isometry.
\end{thm}

Theorem \ref{witt} can easily be obtained by combining Theorem 3.8 on p. 9 with Corollary 5.11 on p.17 of \cite{scharlau}. The first objective of this paper is to make this theorem effective, namely to prove that there exists a decomposition like (\ref{decompose}) with hyperbolic planes, singular, and anisotropic components having relatively small height for an appropriately defined notion of height. By Northcott's theorem, there are only finitely many subspaces of fixed dimension over $K$ whose height is bounded above by a given constant. Hence our result produces a ``search bound'' on components of Witt decomposition for a bilinear space (see \cite{masser:baker} for a discussion of search bounds). This result is also related to the vast collection of results on small-height zeros of quadratic forms. The subject originates in a classical paper of Cassels, \cite{cassels:small}, where he proved that an isotropic rational quadratic form has a zero of relatively small height, producing an explicit bound on height in terms of the height of the quadratic form. Cassels' theorem has been extended and generalized in a number of different ways (see \cite{vaaler:smallzeros}, \cite{masser:baker} and \cite{me:smallzeros} for more information on this). Our first main result can also be viewed in the context of those results; we will discuss this approach in more details in section~3. 
\bigskip

Another direction we pursue here is investigation of the effective structure of the isometry group of a regular symmetric bilinear space $(Z,F)$ over $K$. In \cite{masser:baker} Masser proposes a version of the following question. Let $F$ and $G$ be two symmetric bilinear forms on $K^N$ such that there exists $A \in GL_N(K)$ with $F(A \bX, A \bY) = G(\bX, \bY)$. Can we prove that there exists such an $A$ of bounded height, where the bound would be in terms of heights of $F$ and $G$? In our context $F=G$, and so we can ask for an element of bounded height in the isometry group of the space $(Z,F)$. This question is quite easy to answer (see Corollary~\ref{isom:masser} below), however one can consider the following generalization. Let $\OO(Z,F)$ be the group of isometries of $(Z,F)$. We recall a classical theorem of Cartan and Dieudonn{\'e} (see Theorem 5.4 on p. 15 of \cite{scharlau} or Theorem 43:3 on p.102 of \cite{omeara}). We review the required definitions in section~5.

\begin{thm} \label{CD} Let $(Z,F)$ be a regular symmetric bilinear space over $K$ with $Z \subseteq K^N$ of dimension $L$, $1 \leq L \leq N$. Let $\sigma \in \OO(Z,F)$. Then $\sigma$ can be represented as a product of at most $L$ reflections.
\end{thm}

The identity element of $\OO(Z,F)$ is thought of here as the product of zero reflections. We will be interested in proving a slightly weaker effective version of this theorem, namely given a $\sigma \in \OO(Z,F)$ we will prove that it can be represented as a product of at most $2L-1$ reflections of bounded height, where the bound on height is in terms of heights of $F$, $Z$, and $\sigma$. 
\bigskip   

We start with some notation. We write $d$ for degree of $K$ over $\que$, $O_K$ for its ring of integers, $\D_K$ for its discriminant, and $M(K)$ for its set of places. For each place $v \in M(K)$ we write $K_v$ for the completion of $K$ at $v$ and let $d_v = [K_v:\que_v]$ be the local degree of $K$ at $v$, so that for each $u \in M(\que)$
\begin{equation}
\sum_{v \in M(K), v|u} d_v = d.
\end{equation}

\noindent
For each place $v \in M(K)$ we define the absolute value $\|\ \|_v$ to be the unique absolute value on $K_v$ that extends either the usual absolute value on $\real$ or $\cee$ if $v | \infty$, or the usual $p$-adic absolute value on $\que_p$ if $v|p$, where $p$ is a prime. We also define the second absolute value $|\ |_v$ for each place $v$ by $|a|_v = \|a\|_v^{d_v/d}$ for all $a \in K$. Then for each non-zero $a \in K$ the {\it product formula} reads
\begin{equation}
\label{product_formula}
\prod_{v \in M(K)} |a|_v = 1.
\end{equation} 

\noindent
For each finite place $v \in M(K)$, $v \nmid \infty$, we define the {\it local ring of $v$-adic integers} $O_v = \{ x \in K : |x|_v \leq 1 \}$, whose unique maximal ideal is $P_v =  \{ x \in K : |x|_v < 1 \}$. Then $O_K = \bigcap_{v \nmid \infty} O_v$. For each $v | \infty$ and each positive integer $j$, define as in \cite{vaaler:smallzeros}
\[ r_v(j) = \left\{ \begin{array}{ll}
    \pi^{-1/2} \Gamma(j/2+1)^{1/j} & \mbox{if $v | \infty$ is real} \\
    (2\pi)^{-1/2} \Gamma(j+1)^{1/2j} & \mbox{if  $v | \infty$ is complex}
\end{array}
\right. \]
It will be useful to define a field constant
\begin{equation}
\label{constant:C}
C_K(j) = 2|\D_K|^{1/2d} \prod_{v | \infty} r_v(j)^{d_v/d},
\end{equation}

\noindent
We extend absolute values to vectors by defining the local heights. For each $v \in M(K)$ define a local height $H_v$ on $K_v^N$ by
\[ H_v(\bx) = \left\{ \begin{array}{ll}
\max_{1 \leq i \leq N} |x_i|_v & \mbox{if $v \nmid \infty$} \\
\left( \sum_{i=1}^N \|x_i\|_v^2 \right)^{d_v/2d} & \mbox{if $v | \infty$}
\end{array}
\right. \]
for each $\bx \in K_v^N$. We define the following global height function on $K^N$:
\begin{equation}
H(\bx) = \prod_{v \in M(K)} H_v(\bx),
\end{equation}
for each $\bx \in K^N$. We also define an {\it inhomogeneous} height function on vectors by
 \begin{equation}
h(\bx) = H(1,\bx).
\end{equation}
A basic property of heights that we will use states that for $m_1,...,m_L \in \zed$ and $\bx_1,...,\bx_L \in K^N$,
\begin{equation}
\label{sum_height}
h \left( \sum_{i=1}^L m_i \bx_i \right) \leq \left( \sum_{i=1}^L m_i^2 \right)^{1/2} \prod_{i=1}^L h(\bx_i).
\end{equation}
\smallskip

We extend height to polynomials by viewing it as height function of the coefficient vector of a given polynomial. Hence for our quadratic form $F$, $H(F)$ is the height of the matrix $(f_{ij})_{1 \leq i,j \leq N}$ viewed as a vector in $K^{N^2}$. In general, for an $M \times N$ matrix $A$ we define $H(A)$ by viewing $A$ as a vector in $K^{MN}$, same way as we defined height of~$F$. This way we also have height defined on elements of the isometry group $\OO(K^N,F)$, since they can be represented by $N \times N$ matrices, and each such matrix can be viewed as a vector in $K^{N^2}$. For each element $\sigma$ of the isometry group $\OO(Z,F)$ of a regular bilinear space we will select an extension $\hs \in \OO(K^N,F)$ of minimal possible height, and will define $H(\sigma)$ to be $H(\hs)$. We will explain how this is done in more details in section~5.
\smallskip

We also define another height on matrices, which is the same as height function on subspaces of $K^N$. Let $V \subseteq K^N$ be a subspace of dimension $J$, $1 \leq J \leq N$. Choose a basis $\bx_1,...,\bx_J$ for $V$, and write $X = (\bx_1\ ...\ \bx_J)$ for the corresponding $N \times J$ basis matrix. Then 
$$V = \{ X \bt : \bt \in K^J \}.$$
On the other hand, there exists an $(N-J) \times N$ matrix $A$ with entries in $K$ such that 
$$V = \{ \bx \in K^N : A \bx = 0 \}.$$
Let $\I$ be the collection of all subsets $I$ of $\{1,...,N\}$ of cardinality $J$. For each $I \in \I$ let $I'$ be its complement, i.e. $I' = \{1,...,N\} \setminus I$, and let $\I' = \{ I' : I \in \I\}$. Then 
$$|\I| = \binom{N}{J} = \binom{N}{N-J} = |\I'|.$$
For each $I \in \I$, write $X_I$ for the $J \times J$ submatrix of $X$ consisting of all those rows of $X$ which are indexed by $I$, and $_{I'} A$ for the $(N-J) \times (N-J)$ submatrix of $A$ consisting of all those columns of $A$ which are indexed by $I'$. By the duality principle of Brill-Gordan \cite{gordan:1} (also see Theorem 1 on p. 294 of \cite{hodge:pedoe}), there exists a non-zero constant $\gamma \in K$ such that
\begin{equation}
\label{duality}
\det (X_I) = (-1)^{\varepsilon(I)} \gamma \det (_{I'} A),
\end{equation}
where $\varepsilon(I) = \sum_{i \in I} i$. Define the vectors of {\it Grassmann coordinates} of $X$ and $A$ respectively to be 
$$Gr(X) = (\det (X_I))_{I \in \I} \in K^{|I|},\ \ Gr(A) = (\det (_{I'} A))_{I' \in \I'} \in K^{|I'|}.$$
Define 
$$\H(X) =  H(Gr(X)),\ \ \H(A) = H(Gr(A)),$$
and so by (\ref{duality}) and (\ref{product_formula})
$$\H(X) = \H(A).$$
Define height of $V$ denoted by $H(V)$ to be this common value. Hence the height of a matrix is the height of its row (or column) space, which is equal to the height of its nullspace. Also notice that $Gr(X)$ can be identified with $\bx_1 \wedge\ ...\ \wedge \bx_J$, where $\wedge$ stands for the wedge product, viewed under the cannonical lexicographic embedding into $K^{\binom{N}{J}}$. Therefore we can also write
$$H(V) = H(\bx_1 \wedge\ ...\ \wedge \bx_J).$$
This definition is legitimate, since it does not depend on the choice of the basis for $V$: let $\bwy_1,...,\bwy_J$ be another basis for $V$ over $K$, then there exists $C \in GL_N(K)$ such that $\bwy_i = C \bx_i$ for each $1 \leq i \leq J$, and so
\begin{eqnarray*}
H(\bwy_1 \wedge\ ...\ \wedge \bwy_J) & = & H(C \bx_1 \wedge\ ...\ \wedge C \bx_J) \\
& = & \left( \prod_{v \in M(K)} |\det(C)|_v \right) H(\bx_1 \wedge\ ...\ \wedge \bx_J) \\
& = & H(\bx_1 \wedge\ ...\ \wedge \bx_J),
\end{eqnarray*}
by the product formula. We are now ready to state our main results. First is an effective version of Witt's decomposition Theorem \ref{witt}.

\begin{thm} \label{me_witt} Let $F$ be a symmetric bilinear form on $K^N$. Let $Z \subseteq K^N$ be a subspace of dimension $L$, $2 \leq L \leq N$, and Witt index $M \geq 1$. Let $F$ have rank $r$ on $Z$, $1 \leq r \leq L$. There exists an orthogonal decomposition of the bilinear space $(Z,F)$ of the form (\ref{decompose}) with
\begin{equation}
\label{sing_height}
H(Z^{\perp}) \leq C_K(r)^r H(F)^{r/2} H(Z),
\end{equation}
and
\begin{equation}
\label{height}
\max\{ H(\hyp_i),H(V) \} \leq \A_K(N,L,M) \left\{ H(F)^{\frac{L+2M}{4}} H(Z) \right\}^{\frac{(M+1)(M+2)}{2}},
\end{equation}
for each $1 \leq i \leq M$, where
\begin{equation}
\label{constant}
\A_K(N,L,M) = \left\{ \left( 2^{2M+1} C_K(L)^2 \right)^L \left( N |\D_K|^{1/d} \right)^{M+5L} \right\}^{\frac{M(M+3)}{8}}.
\end{equation}
\end{thm}
\smallskip

Next is an effective version of Cartan-Dieudonn{\'e} Theorem \ref{CD}.

\begin{thm} \label{CD_me} Let $(Z,F)$ be a regular symmetric bilinear space over $K$ with $Z \subseteq K^N$ of dimension $L$, $1 \leq L \leq N$, $N \geq 2$. Let $\sigma \in \OO(Z,F)$. Then either $\sigma$ is the identity, or there exist an integer $1 \leq l \leq 2L-1$ and reflections $\tau_1,...,\tau_l \in \OO(Z,F)$ such that 
\begin{equation}
\label{CD_me_0}
\sigma = \tau_1 \circ \dots \circ \tau_l, 
\end{equation}
and for each $1 \leq i \leq l$,
\begin{equation}
\label{CD_me_1}
H(\tau_i) \leq \left\{ \left( 2 N^2 |\D_K|^{\frac{1}{2d}} \right)^{\frac{L^2}{2}} H(F)^{\frac{L}{3}} H(Z)^{\frac{L}{2}} H(\sigma) \right\}^{5^{L-1}}.
\end{equation}
\end{thm}

\bigskip

This paper is structured as follows. In section~2 we discuss a related problem of producing an orthogonal basis of small height for a bilinear space. This can actually be viewed as a version of Siegel's Lemma for a bilinear space, and provides a decomposition of a bilinear space into an orthogonal sum of one-dimensional subspaces of small height - a result of independent interest. In section~3 we recall some basic lemmas on the properties of bilinear spaces, review a result of Vaaler on a maximal totally isotropic subspace of a bilinear space of small height, and prove an effective decomposition lemma for a bilinear space into a singular and regular components of small height. In section~4 we prove Theorem \ref{me_witt}. In section~5 we develop some notation and preliminary lemmas on the effective structure of the isometry group. In particular, we prove two simple lemmas of independent interest: one on the existence of a small-height isometry of a bilinear space, and the other on the bound for the height of the invariant subspace of an isometry. We use these lemmas in section~6 to prove Theorem \ref{CD_me}. 
\bigskip

\section{Siegel's Lemma for a bilinear space}

In this section we prove a certain analogue of Siegel's Lemma for a bilinear space. First we recall the Bombieri-Vaaler formulation of a general Siegel's Lemma.

\begin{thm}[\cite{vaaler:siegel}] \label{Siegel:gen} Let $U$ be a $J$-dimensional subspace of $K^N$, $J<N$. Then there exists a basis $\bx_1,...,\bx_J \in K^N$ for $U$ such that
\begin{equation}
\label{Siegel:1}
\prod_{i=1}^J H(\bx_i) \leq \prod_{i=1}^J h(\bx_i) \leq \left\{ N |\D_K|^{1/d} \right\}^{J/2} H(U).
\end{equation}
\end{thm}

\noindent
We will also need the following simple technical lemmas.

\begin{lem} \label{intersection} Let $U_1$ and $U_2$ be subspaces of $K^N$. Then
$$H(U_1 \cap U_2) \leq H(U_1) H(U_2).$$
\end{lem}

This well known fact is an immediate corollary of Theorem 1 of \cite{vaaler:struppeck}.

\begin{lem} \label{mult} Let $X$ be a $J \times N$ matrix over $K$ with row vectors $\bx_1,...,\bx_J$, and let $F$ be a symmetric bilinear form in $N$ variables over $K$, as above (we also write $F$ for its $N \times N$ coefficient matrix). Then
$$\H(X F) \leq H(F)^J \prod_{i=1}^J H(\bx_i).$$
\end{lem}

\proof
By Lemma 4.7 of \cite{absolute:siegel}
\begin{equation}
\label{i1}
\H(X F) = H(\bx_1^t F \wedge\ ...\ \wedge \bx_J^t F) \leq \prod_{i=1}^J H(\bx_i^t F).
\end{equation} 
For each $1 \leq i \leq J$,
$$\bx_i^t F = \left( \sum_{j=1}^N f_{j1} x_{ij}, ..., \sum_{j=1}^N f_{jN} x_{ij} \right),$$
and so for $v \nmid \infty$,
\begin{equation}
\label{mx_1}
H_v(\bx_i^t F) \leq H_v(F) H_v(\bx_i),
\end{equation}
and for $v | \infty$, by Cauchy-Schwarz inequality
\begin{eqnarray}
\label{mx_2}
H_v(\bx_i^t F) & = & \left\{ \sum_{k=1}^N \left\| \sum_{j=1}^N f_{jk} x_{ij} \right\|^2 \right\}^{d_v/2d} \nonumber \\
& \leq & \left\{ \sum_{k=1}^N \left( \sum_{j=1}^N \| f_{jk} \|_v^2 \right) \left( \sum_{j=1}^N \| x_{ij} \|_v^2 \right) \right\}^{d_v/2d} = H_v(F) H_v(\bx_i).
\end{eqnarray}
Therefore for each $1 \leq i \leq J$,
\begin{equation}
\label{i2}
H(\bx_i^t F) \leq H(\bx_i) H(F).
\end{equation}
The lemma follows by combining (\ref{i1}) with (\ref{i2}).
\endproof
\bigskip

Next we will use Theorem \ref{Siegel:gen} to produce a small-height orthogonal basis for a subspace of a bilinear space. Specifically, we prove the following theorem.

\begin{thm} \label{Siegel:me} Let $U$ be a $J$-dimensional subspace of $(K^N,F)$, $J<N$. Then there exists a basis $\bx_1,...,\bx_J \in K^N$ for $U$ such that $F(\bx_i,\bx_j) = 0$ for all $i \neq j$, and
\begin{equation}
\label{Siegel:2}
\prod_{i=1}^J H(\bx_i) \leq \left( N |\D_K| \right)^{\frac{J^2+J-2}{4}} H(F)^{\frac{J(J+1)}{2}} H(U)^J.
\end{equation}
\end{thm}

\proof
We argue by induction on $J$. First suppose that $J=1$, then pick any $\boldsymbol 0 \neq \bx_1 \in U$, and observe that $H(\bx_1) = H(U)$. Now assume that $J>1$ and the theorem is true for all $1 \leq j < J$. Let $\boldsymbol 0 \neq \bx_1 \in U$ be a vector guaranteed by Theorem \ref{Siegel:gen} so that
\begin{equation}
\label{s2}
H(\bx_1) \leq \left\{ N |\D_K|^{1/d} \right\}^{1/2} H(U)^{1/J}.
\end{equation}
First assume that $\bx_1$ is a non-singular point in $U$. Then
$$U_1 = \{ \bwy \in U : \bx_1^t F \bwy = 0 \} = \{ \bx_1 \}^{\perp} \cap U,$$
has dimension $J-1$; here $\{ \bx_1 \}^{\perp} = \{ \bwy \in K^N : \bx_1^t F \bwy = 0 \}$. Then by Lemma \ref{intersection}, Lemma \ref{mult}, and (\ref{s2}) we obtain
\begin{equation}
\label{s3}
H(U_1) \leq H(\bx_1^t F) H(U) \leq H(F) H(\bx_1) H(U) \leq \left( N |\D_K|^{1/d} \right)^{1/2} H(F) H(U)^{\frac{J+1}{J}}.
\end{equation}
Since $\dim_K(U_1) = J-1$, the induction hypothesis implies that there exists a basis $\bx_2,...,\bx_{J}$ for $U_1$ such that $F(\bx_i,\bx_j) = 0$ for all $2 \leq i \neq j \leq J$, and
\begin{eqnarray}
\label{s4} 
\prod_{i=2}^J H(\bx_i) & \leq & \left( N |\D_K|^{1/d} \right)^{\frac{J^2-J-2}{4}} H(F)^{\frac{J(J-1)}{2}} H(U_1)^{J-1} \nonumber \\
& \leq & \left( N |\D_K|^{1/d} \right)^{\frac{J^2+J-4}{4}} H(F)^{\frac{J^2+J-2}{2}} H(U)^{\frac{J^2-1}{J}},
\end{eqnarray}
where the last inequality follows by (\ref{s3}). Combining (\ref{s2}) and (\ref{s4}) we see that $\bx_1,...,\bx_J$ is a basis for $U$ satisfying (\ref{Siegel:2}) such that $F(\bx_i,\bx_j) = 0$ for all $1 \leq i \neq j \leq J$.

Now assume that $\bx_1$ is a singular point in $U$. Since $\bx_1 \neq 0$, it must be true that $x_{1j} \neq 0$ for some $1 \leq j \leq N$. Let
$$U_1 = U \cap \{ \bx \in K^N : x_j = 0 \},$$
then $\bx_1 \notin U_1$, $U = K\bx_1 \perp U_1$, and 
\begin{equation}
\label{s5}
H(U_1) \leq H(U),
\end{equation}
by Lemma \ref{intersection}. Since $\dim_K(U_1) = J-1$, we can apply induction hypothesis to $U_1$, and proceed the same way as in the non-singular case above. Since the upper bound of (\ref{s5}) is smaller than that of (\ref{s3}), the result follows.
\endproof
\smallskip

Notice that Theorem \ref{Siegel:me} can be reformulated by saying that there exists a decomposition of the bilinear space $(U,F)$ into an orthogonal sum of one-dimensional subspaces, the product of heights of which is bounded above by (\ref{Siegel:2}). Therefore Theorem \ref{Siegel:me} can also be viewed as a result on effective orthogonal decomposition of a bilinear space, which is the subject of this paper.  
\bigskip

\section{Small zeros of quadratic forms}

Let $F$ be a symmetric bilinear form in $2N$ variables over $K$, as above. Let $Z \subseteq K^N$ be a subspace of dimension $2 \leq L \leq N$. We write $(Z,F)$ for the bilinear space on $Z$ with the bilinear form $F$ restricted to $Z$. In this section we review some basic results on bilinear spaces and setup the notation that will later be used in the proof of Theorem \ref{me_witt}.
\smallskip

We start by giving a brief overview of required notation (see Chapter 1 of \cite{scharlau} for a detailed introduction into the subject). A totally isotropic subspace $W$ of $(Z,F)$ is a subspace such that for all $\bx,\bwy \in W$, $F(\bx,\bwy)=0$. All maximal totally isotropic subspaces of $(Z,F)$ have the same dimension. It is called the Witt index of $(Z,F)$ and we denote it by $M$. A subspace $U$ of $(Z,F)$ is anisotropic if $F(\bx) \neq 0$ for all $\boldsymbol 0 \neq \bx \in U$. A subspace $U$ of $(Z,F)$ is called regular if for each $\boldsymbol 0 \neq \bx \in U$ there exists $\bwy \in U$ so that $F(\bx,\bwy) \neq 0$. For each subspace $U$ of $(Z,F)$ we define $U^{\perp} = \{ \bx \in Z : F(\bx, \bwy) = 0\ \forall\ \bwy \in U \}$. If two subspaces $U_1$ and $U_2$ of $(Z,F)$ are orthogonal, we write $U_1 \perp U_2$ for their orthogonal sum. If $U$ is a regular subspace of $(Z,F)$, then $Z = U \perp U^{\perp}$ and $U \cap U^{\perp} = \{\boldsymbol 0\}$.
\smallskip

Two vectors $\bx,\bwy \in Z$ are called a hyperbolic pair if $F(\bx) = F(\bwy) = 0$, $F(\bx,\bwy) = 1$; the subspace $\hyp(\bx,\bwy) = \spn_K \{\bx,\bwy\}$ is regular and is called a hyperbolic plane. An orthogonal sum of hyperbolic planes is called a hyperbolic space. Every hyperbolic space is regular.
\smallskip

We now state a result of Vaaler \cite{vaaler:smallzeros} (see also \cite{schmidt:schlickewei}) on the existence of a maximal totally isotropic subspace of $(Z,F)$ of small height, which we later use in the proof of Theorem \ref{me_witt}.

\begin{thm} [\cite{vaaler:smallzeros}] \label{vaaler:quadratic} Let $M \geq 1$ be the Witt index of $(Z,F)$ over $K$. Then there exists a subspace $W$ of $(Z,F)$ of dimension $M$ such that $F(\bx)=0$ for all $\bx \in W$ and
\begin{equation}
\label{vaaler1}
H(W) \leq \{2^{2M+1} C_K(L-M)^2 H(F)\}^{(L-M)/2} H(Z).
\end{equation}  
\end{thm}

Notice that subspace $W$ of Theorem \ref{vaaler:quadratic} is indeed maximal totally isotropic. Maximality is by construction. Also, for each $\bx,\bwy \in W$, $\bx + \bwy \in W$, hence
$$0 = F(\bx+\bwy) = F(\bx) + F(\bwy) + 2F(\bx,\bwy) = 2F(\bx,\bwy).$$
\smallskip

A consequence of a related theorem of Vaaler is the following simple decomposition lemma in case when $(Z,F)$ is not a regular space.

\begin{lem} \label{reg:decompose} Let $F$ have rank $r$ on $Z$, and assume that $1 \leq r < L$. Then the bilinear space $(Z,F)$ can be represented as
\begin{equation}
\label{regular1}
Z = Z^{\perp} \perp W,
\end{equation}
where $W$ is a regular subspace of $Z$, with
\begin{equation}
\label{regular2}
H(Z^{\perp}) \leq C_K(r)^r H(F)^{r/2} H(Z),
\end{equation}
and
\begin{equation}
\label{regular3}
H(W) \leq \left\{ N |\D_K|^{1/d} \right\}^{L/2} H(Z).
\end{equation}
\end{lem}

\proof
The fact that $Z^{\perp}$ satisfies (\ref{regular2}) is guaranteed by Theorem 2 of \cite{vaaler:smallzeros2}. Now let $\bz_1,...,\bz_L$ be the basis for $Z$ guranteed by Theorem \ref{Siegel:gen}, then
\begin{equation}
\label{r1}
\prod_{i=1}^L H(\bz_i) \leq \left\{ N |\D_K|^{1/d} \right\}^{L/2} H(Z).
\end{equation}
Notice that $\dim_K(Z^{\perp})=L-r$. We can now pick $r$ vectors $\bz_{i_1},...,\bz_{i_r}$ from our basis for $Z$ such that
$$\spn_K \{Z^{\perp},\bz_{i_1},...,\bz_{i_r}\} = Z.$$
Let $W = \spn_K \{\bz_{i_1},...,\bz_{i_r}\}$, then $Z = Z^{\perp} \oplus W$. This implies, by Theorem 3.8 on p. 9 of \cite{scharlau}, that $Z = Z^{\perp} \perp W$, $W$ is regular and unique up to isometry. Also, combining Lemma 4.7 of \cite{absolute:siegel} with (\ref{r1}), we obtain
$$H(W) = H(\bz_{i_1} \wedge\ ...\ \wedge \bz_{i_{r}}) \leq \prod_{j=1}^r H(\bz_{i_j}) \leq \left\{ N |\D_K|^{1/d} \right\}^{L/2} H(Z).$$
This finishes the proof.
\endproof
\smallskip

Notice that we can immediately deduce a version of Cassels' theorem on small zeros of quadratic form $F$ over $K$ from Theorem \ref{vaaler:quadratic}. Namely, if $F$ is isotropic over $K$, then there exists $\boldsymbol 0 \neq \bx \in \V_K(F) = \{ \bt \in K^N : F(\bt) = 0 \}$ such that
\begin{equation}
\label{q1}
H(\bx) \ll_{K,N} H(F)^{\frac{N-1}{2}}.
\end{equation}
The exponent $\frac{N-1}{2}$ on $H(F)$ is proved to be best possible. In fact, if $\V_K(F)$ contains a {\it nonsingular} point, then by Corollary 1.2 of \cite{me:smallzeros} there exists such a point satisfying (\ref{q1}). A similar statement about singular points of small height in $\V_K(F)$ can be deduced from Lemma \ref{reg:decompose}.

\begin{cor} \label{singular} Suppose that $\V_K(F) = \{ \bt \in K^N : F(\bt) = 0 \}$ contains a singular point $\bx \neq \boldsymbol 0$, so $1 \leq r = \rk(F) < N$. Then there exists such a point $\bx$ with
\begin{equation}
\label{sing}
H(\bx) \leq \sqrt{N} |\D_K|^{1/2d} C_K(r)^{\frac{r}{N-r}} H(F)^{\frac{r}{2(N-r)}}.
\end{equation} 
\end{cor}

\proof
Let $Z$ of Lemma \ref{reg:decompose} be $K^N$, then $H(Z)=1$, $L=N$, and $\dim_K(Z^{\perp})=N-r$. Clearly $Z^{\perp} \subseteq \V_K(F)$, and all points of $Z^{\perp}$ are singular in $\V_K(F)$. By Theorem \ref{Siegel:gen}, there must exist $\boldsymbol 0 \neq \bx \in Z^{\perp}$ such that
$$H(\bx) \leq \sqrt{N} |\D_K|^{1/2d} H(Z^{\perp})^{1/(N-r)} \leq \sqrt{N} |\D_K|^{1/2d} C_K(r)^{\frac{r}{N-r}} H(F)^{\frac{r}{2(N-r)}},$$
where the last inequality follows by (\ref{regular2}).
\endproof

Notice that Corollary \ref{singular} suggests that in this context the singular case can be simpler than the nonsingular one. This unusual phenomenon has already been observed in \cite{masser:1} and \cite{me:smallzeros}. We are now ready to prove Theorem \ref{me_witt}.
\bigskip

\section{Proof of Theorem \ref{me_witt}}

We first prove a version of our theorem for a regular bilinear space. We remark that everywhere in our arguments, if $m < n$, then $\sum_{i=n}^m$ is taken to mean $0$ and $\prod_{i=n}^m$ is taken to mean $1$.

\begin{thm} \label{me_witt_reg} Let $F$ be a symmetric bilinear form on $K^N$. Let $Z \subseteq K^N$ be a subspace of dimension $L$, $2 \leq L \leq N$, such that the bilinear space $(Z,F)$ is regular, i.e. $Z^{\perp}=\{ \boldsymbol 0 \}$. Let $M \geq 1$ be the Witt index of $(Z,F)$. There exists an orthogonal decomposition of $(Z,F)$ of the form 
\begin{equation}
\label{reg_dec}
Z = \hyp_1 \perp\ ...\ \perp \hyp_M \perp V,
\end{equation}
where $\hyp_i$ are hyperbolic planes, $V$ is anisotropic component, and
\begin{equation}
\label{height_reg}
\max\{ H(\hyp_i),H(V) \} \leq A_K(N,L,M) \left\{ H(F)^{\frac{L+2M}{4}} H(Z) \right\}^{\frac{(M+1)(M+2)}{2}},
\end{equation}
for each $1 \leq i \leq M$, where
\begin{equation}
\label{constant_reg}
A_K(N,L,M) = \left\{ \left( 2^{2M+1} C_K(L)^2 \right)^L \left( N |\D_K|^{1/d} \right)^{M+L} \right\}^{\frac{M(M+3)}{8}}.
\end{equation}
\end{thm}

\proof
Let $W$ be a maximal totally isotropic subspace of $(Z,F)$ satisfying (\ref{vaaler1}) and let $\bx_1,...,\bx_M$ be the basis for $W$ guaranteed by Theorem \ref{Siegel:gen}. Notice that $F(\bx_i,\bx_j) = 0$ for all $1 \leq i,j \leq M$, since $W$ is a totally isotropic subspace. Let $\bwy_1,...,\bwy_L$ be the basis for $Z$ guaranteed by Theorem \ref{Siegel:gen}, ordered so that 
$$H(\bwy_1) \leq H(\bwy_2) \leq ... \leq H(\bwy_L).$$
For each $1 \leq i \leq M$ let $j_i$ be the smallest index such that $F(\bx_i,\bwy_{j_i}) \neq 0$. Such $j_i$ exists for each $i$ since otherwise $\bx_i$ would be a singular point, contradicting regularity of $(Z,F)$. By reordering $\bx_1,...,\bx_M$ if necessary, we can assume without loss of generality that
$$1 \leq j_M \leq j_{M-1} \leq ... \leq j_1 \leq L.$$
Moreover, for each $1 \leq i \leq M$, $j_i \leq L-i+1$, since
$$\spn_K \{ \bwy_1,...,\bwy_{L-i+1} \} \nsubseteq \spn_K \{ \bx_1,...,\bx_i \}^{\perp},$$   
and so $H(\bwy_{j_i}) \leq H(\bwy_{L-i+1})$ by our ordering of $\bwy_1,...,\bwy_L$. Therefore, by (\ref{Siegel:1})
\begin{eqnarray}
\label{pr1}
\prod_{i=1}^M H(\bx_i) H(\bwy_{j_i}) & \leq & \prod_{i=1}^M H(\bx_i) H(\bwy_{L-i+1}) \nonumber \\
& = & \left( \prod_{i=1}^M H(\bx_i) \right) \left( \prod_{i=1}^M H(\bwy_{L-i+1}) \right) \nonumber \\
& \leq & \left\{ N |\D_K|^{1/d} \right\}^{\frac{M+L}{2}} H(W) H(Z).
\end{eqnarray} 
In particular, for some $1 \leq i \leq M$, we must have
\begin{equation}
\label{pr2}
H(\bx_i) H(\bwy_{j_i}) \leq \left\{ N |\D_K|^{1/d} \right\}^{\frac{M+L}{2M}} \left( H(W) H(Z) \right)^{\frac{1}{M}}.
\end{equation}
Define $\hyp_1 = \spn_K \{ \bx_i,\bwy_{j_i} \}$ for this choice of $i$. Since $F(\bx_i)=0$ and $F(\bx_i,\bwy_{j_i}) \neq 0$, $\hyp_1$ is a regular subspace of $Z$ with Witt index equal to one, hence it is a hyperbolic plane. Notice that by combining (\ref{pr2}) and (\ref{vaaler1}), we have
\begin{equation}
\label{pr3}
H(\hyp_1) \leq H(\bx_i) H(\bwy_{j_i}) \leq  B_K(N,L,M) H(F)^{\frac{L-M}{2M}} H(Z)^{\frac{2}{M}},
\end{equation}
where
\begin{equation}
\label{const}
B_K(N,L,M) = \left\{ \left( 2^{2M+1} C_K(L-M)^2 \right)^{L-M} \left( N |\D_K|^{1/d} \right)^{M+L} \right\}^{\frac{1}{2M}}.
\end{equation}
Define 
$$Z_1 = \hyp_1^{\perp} = \{ \bz \in K^N : F(\bz,\bx) = 0\ \forall\ \bx \in \hyp_1 \} \cap Z,$$
so $\dim_K(Z_1) = L-2$, and $Z = \hyp_1 \perp Z_1$. Notice that by combining Lemma \ref{intersection}, Lemma \ref{mult}, and (\ref{pr3}), we have
\begin{equation}
\label{pr4}
H(Z_1) \leq H(\hyp_1) H(Z) H(F)^2 \leq B_K(N,L,M) H(F)^{\frac{L+3M}{2M}} H(Z)^{\frac{M+2}{M}}.
\end{equation}
We continue by induction on $M$. If $M=1$, we are done. If $M \geq 2$, assume that the theorem holds for a bilinear space of Witt index smaller than $M$, in particular it holds for $(Z_1,F)$, a bilinear space of dimension $L-2$ and Witt index $M-1$. Then there exists a decomposition
\begin{equation}
\label{pr5}
Z_1 = \hyp_2 \perp\ ...\ \perp \hyp_M \perp V,
\end{equation}
where $V$, the anisotropic component of $Z_1$ is the same as that of $Z$, and combining the induction hypothesis with (\ref{pr4}) and (\ref{const}), for each $2 \leq i \leq M$ we obtain
\begin{eqnarray}
\label{pr6}
\max\{ H(\hyp_i),H(V) \} & \leq & A_K(N,L-2,M-1) \left\{ H(F)^{\frac{L+2M-4}{4}} H(Z_1) \right\}^{\frac{M(M+1)}{2}} \nonumber \\
& \leq & A_K(N,L-2,M-1) B_K(N,L,M)^{\frac{M(M+1)}{2}} \times \nonumber \\
& \times & \left\{ H(F)^{\frac{L+2M-4}{4} + \frac{L+3M}{2M}} H(Z)^{\frac{M+2}{M}} \right\}^{\frac{M(M+1)}{2}} \nonumber \\
& \leq & A_K(N,L,M) \left\{ H(F)^{\frac{L+2M}{4}} H(Z) \right\}^{\frac{(M+1)(M+2)}{2}}.
\end{eqnarray}
This finishes the proof.
\endproof
\smallskip

\noindent
{\it Proof of Theorem \ref{me_witt}.} If $(Z,F)$ is regular, then $Z^{\perp} = \{\boldsymbol 0\}$, and we are done by Theorem \ref{me_witt_reg}. Let $r$ be rank of $F$ on $Z$, and assume that $1 \leq r < L$. By Lemma \ref{reg:decompose}, there exists a decomposition of $Z$ of the form (\ref{regular1}) with $H(Z^{\perp})$ and $H(W)$ bounded as in (\ref{regular2}) and (\ref{regular3}) respectively. Now $W$ is a regular subspace of $Z$, so we can apply Theorem \ref{me_witt_reg} to the bilinear space $(W,F)$. The result follows.
\boxed{ }
\bigskip

\section{Isometries of a bilinear space}

In this section we develop the preliminaries needed for the proof of Theorem \ref{CD_me}. We start with some definitions and then prove a few technical lemmas. Let $F$ be a symmetric bilinear form as above, and let $Z$ be an $L$-dimensional subspace of $K^N$, $1 \leq L \leq N$, $N \geq 2$, such that the bilinear space $(Z,F)$ is regular, and thus $K^N = Z \perp Z^{\perp_{K^N}}$, where $Z^{\perp_{K^N}} = \{ \bx \in K^N : F(\bx,\bz) = 0\ \forall\ \bz \in Z\}$. Let $\OO(Z,F)$ be the group of isometries of $(Z,F)$, and write $id_Z$ for its identity element. Also let $-id_Z$ be the element of $\OO(Z,F)$ that takes $\bx$ to $-\bx$ for each $\bx \in Z$. Each element $\sigma$ of the isometry group $\OO(K^N,F)$ is uniquely represented by an $N \times N$ matrix $A \in GL_N(K)$, and so we can define $H(\sigma) = H(A)$, where $H(A)$ is defined by viewing $A$ as a vector in $K^{N^2}$ as we did in section~1. 

Notice that each $\sigma \in \OO(Z,F)$ can be extended to an isometry of $\hs \in \OO(K^N,F)$ by selecting an isometry $\sigma' \in \OO(Z^{\perp_{K^N}},F)$. For each $\sigma \in \OO(Z,F)$ choose such an extension $\hs: K^N \rightarrow K^N$ so that $H(\hs)$ is minimal, and define $H(\sigma) = H(\hs)$ for this choice of $\hs$. This definition of height in particular insures that for each $\sigma \in \OO(Z,F)$
\begin{equation}
\label{pm}
H(\sigma) = H(-\sigma),
\end{equation}
where $-\sigma = -id_Z \circ \sigma$. Moreover, if $A$ is the matrix of $\hs$, then
\begin{equation}
\label{det}
\det(A) = \det (\hs) = \det \left( \hs \mid_Z \right) \det \left( \hs \mid_{Z^{\perp_{K^N}}} \right) = \det(\sigma) \det(\sigma') = \pm 1.
\end{equation}
We will also refer to this matrix $A$ as the matrix of $\sigma$.

For each $\bx \in Z$ such that $F(\bx) \neq 0$ we can define an element of $\OO(Z,F)$, $\tau_{\bx} : Z \longrightarrow Z$, given by
\begin{equation}
\label{cd1}
\tau_{\bx} (\bwy) = \bwy - \frac{2 F(\bx,\bwy)}{F(\bx)} \bx,
\end{equation}
which is a {\it reflection} in the hyperplane $\{\bx\}^{\perp} = \{ \bz \in Z : F(\bx,\bz) = 0 \}$. It is not difficult to see that the matrix of such a reflection is of the form $(\tau_{ij}(\bx))_{1 \leq i,j \leq N}$, where
\[ \tau_{ij}(\bx) = \left\{ \begin{array}{ll}
    1 - \frac{2}{F(\bx)} \sum_{k=1}^N f_{ik} x_i x_k & \mbox{if $i=j$} \\
    - \frac{2}{F(\bx)} \sum_{k=1}^N f_{jk} x_i x_k & \mbox{if $i \neq j$}
    \end{array}
\right. \]
For each reflection $\tau_{\bx}$, $\det(\tau_{\bx}) = -1$. We say that $\sigma$ is a {\it rotation} if $\det(\sigma)=+1$.

\begin{lem} \label{ref_height} Let $\bx \in Z$ be anisotropic and $\tau_{\bx} \in \OO(Z,F)$ be the corresponding reflection. Then
\begin{equation}
\label{ref_bound}
H(\tau_{\bx}) \leq N^3(N+2) H(F) H(\bx)^2.
\end{equation}
\end{lem}

\proof
By the product formula, $H(\tau_{\bx}) = H(F(\bx) \tau_{\bx})$. If $v \in M(K)$ is such that $v \nmid \infty$, then for each $1 \leq i=j \leq N$
\begin{eqnarray}
\label{cd2}
|F(\bx) \tau_{ij}(\bx)|_v & = & \left| F(\bx) - 2 \sum_{k=1}^N f_{jk} x_i x_k \right|_v \nonumber \\
& = & \left| \sum_{l=1}^N \sum_{m=1}^N f_{lm} x_l x_m - 2 \sum_{k=1}^N f_{jk} x_i x_k \right|_v \leq H_v(F) H_v(\bx)^2,
\end{eqnarray}
since $|2|_v \leq 1$, and similarly when $i \neq j$, so $H_v(F(\bx) \tau_{\bx}) \leq H_v(F) H_v(\bx)^2$.
\smallskip

\noindent
If $v | \infty$, then for each $1 \leq i=j \leq N$
\begin{eqnarray}
\label{cd3}
\|F(\bx) \tau_{ij}(\bx)\|_v & \leq & \sum_{l=1}^N \sum_{m=1}^N \|f_{lm} x_l x_m\|_v + 2 \sum_{k=1}^N \|f_{jk} x_i x_k\|_v \nonumber \\
& \leq & N(N+2) \max_{1 \leq l,m \leq N} \|f_{lm} x_l x_m\|_v \nonumber \\
& \leq & N(N+2) \left\{ H_v(F) H_v(\bx)^2 \right\}^{d/d_v},
\end{eqnarray}
and similarly when $i \neq j$, therefore $H_v(F(\bx) \tau_{\bx}) \leq \left\{ N^3(N+2) \right\}^{d_v/d} H_v(F) H_v(\bx)^2$. The result follows by taking a product over all places of $K$.
\endproof

\begin{lem} \label{nonzero} Let $\sigma \in \OO(Z,F)$. There exists an anisotropic vector $\bwy$ in $Z$ such that $\sigma(\bwy) \pm \bwy$ is also anisotropic for some choice of $\pm$, and
\begin{equation}
\label{anis}
H(\bwy) \leq h(\bwy) \leq 2 \sqrt{L} \left\{ N |\D_K|^{\frac{1}{d}} \right\}^{\frac{L+2}{4}} H(Z)^{\frac{L+2}{2L}}.
\end{equation}
\end{lem}

\proof
If $L=1$, then $Z=K\bwy$ for some $\boldsymbol 0 \neq \bwy \in K^N$, and since $(Z,F)$ is regular, $F(\bwy) \neq 0$, $H(\bwy) = H(Z)$, $\OO(Z,F)=\{ id_Z \}$, and clearly $id_Z(\bwy)+\bwy = 2\bwy$ is also anisotropic. Hence assume $L \geq 2$. Let $\bx_1,...,\bx_L$ be a basis for $Z$ which satisfies (\ref{Siegel:1}), ordered so that
$$h(\bx_1) \leq h(\bx_2) \leq \dots \leq h(\bx_L),$$
Let $m$ be the the smallest index such that the restriction of $F$ to 
$$U = \spn_K \{\bx_1,...,\bx_m\}$$
is not identically zero. Since $(Z,F)$ is regular, we must have $1 \leq m \leq \left[ \frac{L}{2} \right] + 1$, and therefore, by (\ref{Siegel:1})
\begin{equation}
\label{ord}
\prod_{i=1}^m h(\bx_i) \leq \left\{ N |\D_K|^{\frac{1}{d}} \right\}^{\frac{m}{2}} H(Z)^{\frac{m}{L}} \leq \left\{ N |\D_K|^{\frac{1}{d}} \right\}^{\frac{L+2}{4}} H(Z)^{\frac{L+2}{2L}}.
\end{equation} 
Notice that for every vector $\bx \in Z$,
$$F(\sigma(\bx)-\bx) + F(\sigma(\bx)+\bx) = 4F(\bx).$$
Since $F$ is not identically zero on $U$, it must therefore be true that at least one of $F \circ (\sigma \pm id_Z)$ is not identically zero on $U$. Assume for instance that $F \circ (\sigma - id_Z)$ is not identically zero on $U$. Then the homogeneous polynomial of degree four in $m$ variables 
$$P(a_1,...,a_m) = F \left( \sum_{i=1}^m a_i \bx_i \right) F \left( \sigma \left( \sum_{i=1}^m a_i \bx_i \right) - \sum_{i=1}^m a_i \bx_i\right) \in K[a_1,...,a_m]$$
is not identically zero on $U$. Therefore there exist $\beta_1,...,\beta_m \in \{-2,-1,0,1,2\}$ such that $P(\beta_1,...,\beta_m) \neq 0$. Let $\bwy = \sum_{i=1}^m \beta_i \bx_i$ for this choice of $\beta_1,...,\beta_m$, then $\bwy \in U$ is precisely the vector we are looking for. Combining (\ref{sum_height}) and (\ref{ord}) we obtain
$$H(\bwy) \leq h(\bwy) \leq \sqrt{\frac{4(L+2)}{2}} \prod_{i=1}^m h(\bx_i) \leq 2 \sqrt{L} \left\{ N |\D_K|^{\frac{1}{d}} \right\}^{\frac{L+2}{4}} H(Z)^{\frac{L+2}{2L}},$$
since $L \geq 2$. This completes the proof.
\endproof
\smallskip

An immediate consequence of Lemma \ref{ref_height} and Lemma \ref{nonzero} is the following statement on the existence of isometries of $(Z,F)$ of small height. This is related to a question of Masser in \cite{masser:baker} (see the discussion on this in the introduction - section~1).

\begin{cor} \label{isom:masser} There exists a reflection $\tau \in \OO(Z,F)$ with
\begin{equation}
\label{isom_bound}
H(\tau) \leq 4L N^{\frac{L+8}{2}} (N+2) |\D_K|^{\frac{L+2}{2d}} H(F) H(Z)^{\frac{L+2}{L}}.
\end{equation}
\end{cor}

\proof
Let $\bx$ be an anisotropic point in $Z$ guaranteed by Lemma \ref{nonzero}. Let $\tau=\tau_{\bx}$. The result follows by combining (\ref{ref_bound}) with (\ref{anis}).
\endproof

\begin{lem} \label{matrix} Let $A \in GL_N(K)$ be such that $\det(A) = \pm 1$, and write $I_N$ for the $N \times N$ identity matrix. Then
\begin{equation}
\label{matrix_height}
H(A \pm I_N) \leq 2 H(A).
\end{equation}
\end{lem}

\proof
Let $\baa_1,...,\baa_N$ be row vectors of $A$. Then for each $v \in M(K)$
\[ \prod_{i=1}^N H_v(\baa_i) \geq \left\{ \begin{array}{ll}
     | \det(A) |_v = 1 & \mbox{if $v \nmid \infty$} \\
     \| \det(A) \|_v = 1 & \mbox{if  $v | \infty$}
\end{array}
\right. \]
by Hadamard's inequality. Therefore, if $v \nmid \infty$, we have 
$$H_v(A) = \max_{1 \leq i \leq N} \{H_v(\baa_i)\} \geq 1,$$ 
and so 
\begin{equation}
\label{m1}
H_v(A \pm I_N) \leq \max \{ 1,H_v(A) \} = H_v(A).
\end{equation}
If $v | \infty$,
\begin{eqnarray*}
1 \leq \left( \prod_{i=1}^N H_v(\baa_i)^{\frac{d}{d_v}} \right)^{\frac{1}{N}} & \leq & \frac{1}{N} \sum_{i=1}^N H_v(\baa_i)^{\frac{d}{d_v}} \\
& \leq & \frac{1}{\sqrt{N}} \left( \sum_{i=1}^N H_v(\baa_i)^{\frac{2d}{d_v}} \right)^{\frac{1}{2}} =  \frac{1}{\sqrt{N}} H_v(A)^{\frac{d}{d_v}},
\end{eqnarray*}
where the last inequality follows by Cauchy-Schwarz. Hence $H_v(A)^{\frac{d}{d_v}} \geq \sqrt{N}$, and so, by the triangle inequality,
\begin{equation}
\label{m2}
H_v(A \pm I_N)^{\frac{d}{d_v}} \leq H_v(A)^{\frac{d}{d_v}} + H_v(I_N)^{\frac{d}{d_v}} \leq H_v(A)^{\frac{d}{d_v}} + \sqrt{N} \leq 2 H_v(A)^{\frac{d}{d_v}}.
\end{equation}
The result follows by combining (\ref{m1}) with (\ref{m2}) and taking a product over all places of $K$.
\endproof
\smallskip

The following simple corollary of Lemma \ref{matrix} provides a bound on the height of the invariant subspace of an isometry, which is an object of interest in the algebraic theory of quadratic forms.

\begin{cor} \label{invariant} Let $\sigma \in \OO(Z,F)$. Let $U$ be the invariant subspace of $\sigma$, i.e. $U = \{ \bz \in Z : \sigma(\bz) = \bz \}$. Let $J = \dim_K(U) \leq L$. Then
\begin{equation}
\label{inv_height}
H(U) \leq \left\{ 2 H(\sigma) \right\}^{N-J} H(Z),
\end{equation}
\end{cor}

\proof
Write $A$ for the $N \times N$ matrix of $\sigma$, and $I_N$ for the $N \times N$ identity matrix. Notice that $U = \{ \bz \in Z : (A-I_N) \bz = \boldsymbol 0 \}$. Let $B$ be a submatrix of $A-I_N$ which consists of $N-J$ linearly independent rows of $A-I_N$. Hence rows of $B$ are of the form $\baa_{i_1}-\be_{i_1},...,\baa_{i_{N-J}}-\be_{i_{N-J}}$ for some $i_1,...,i_{N-J} \in \{1,...,N\}$. Then, by Lemma 4.7 of \cite{absolute:siegel}
\begin{eqnarray}
\label{cd5}
\H(B) & = & H\left( (\baa_{i_1}-\be_{i_1}) \wedge\ ...\ \wedge (\baa_{i_{N-J}}-\be_{i_{N-J}}) \right) \nonumber \\
& \leq & \prod_{j=1}^{N-J} H(\baa_{i_j}-\be_{i_j}) \leq H(A-I_N)^{N-J} \leq (2 H(A))^{N-J},
\end{eqnarray}
where the last inequality follows by Lemma \ref{matrix}. Combining (\ref{cd5}) with Lemma \ref{intersection}, we obtain
$$H(U) \leq \H(B) H(Z) \leq \{ 2 H(A) \}^{N-J} H(Z).$$ 
The finishes the proof, since $H(\sigma) = H(A)$ by definition.
\endproof
\smallskip

The following lemma bounds the height of a product of two matrices.

\begin{lem} \label{matrix1} Let $A$ and $B$ be two $N \times N$ matrices with entries in $K$. Then
\begin{equation}
\label{mht}
H(AB) \leq H(A) H(B).
\end{equation}
\end{lem}

\proof
Write $A = ( \baa_1 \dots \baa_N)^t$, i. e. $\baa_1^t, \dots, \baa_N^t$ are row vectors of $A$. Then we can think of
$$AB = (\baa_1^t B, \dots, \baa_N^t B)^t$$
as a vector in $K^{N^2}$. Hence for each $v \in M(K)$ such that $v \nmid \infty$
$$H_v(AB) = \max_{1 \leq i \leq N} \{ H_v(\baa_i^t B) \} \leq H_v(B) \max_{1 \leq i \leq N} \{ H_v(\baa_i) \} = H_v(A) H_v(B),$$
by (\ref{mx_1}). For each $v | \infty$, we have
$$H_v(AB) = \left\{ \sum_{i=1}^N H_v(\baa_i^t B)^{\frac{2d}{d_v}} \right\}^{\frac{d_v}{2d}} \leq H_v(B) \left\{ \sum_{i=1}^N H_v(\baa_i)^{\frac{2d}{d_v}} \right\}^{\frac{d_v}{2d}} = H_v(A) H_v(B),$$
by (\ref{mx_2}). The conclusion follows by taking a product.
\endproof

\bigskip

\section{Effective version of Cartan-Dieudonn{\'e} theorem}

In this section we will prove Theorem \ref{CD_me}. Let all the notation be as in section~5. We argue by induction on $L$. When $L=1$, $Z=K\bx$ for some anisotropic vector $\bx \in K^N$, since $(Z,F)$ is regular. Then $\sigma = \pm id_Z$, where $-id_Z = \tau_{\bx}$, and $H(\sigma) = \sqrt{N}$ by (\ref{pm}). 
\smallskip

Then assume $L>1$. Write $A$ for the $N \times N$ matrix of $\sigma$, and $I_N$ for the $N \times N$ identity matrix, so in particular $H(\sigma)=H(A)$. Notice that for each $\bx \in Z$,
\begin{equation}
\label{zero}
F(\sigma(\bx)-\bx,\sigma(\bx)+\bx) = 0.
\end{equation}
Let $\bx \in Z$ be the anisotropic vector guaranteed by Lemma \ref{nonzero} with $\sigma(\bx) \pm \bx$ also anisotropic. For this choice of $\pm$, $\tau_{\sigma(\bx) \pm \bx}$ fixes $\sigma(\bx) \mp \bx$ and maps $\sigma(\bx) \pm \bx$ to $-(\sigma(\bx) \pm \bx)$. Then $2 \sigma(\bx) = (\sigma(\bx)) + (\sigma(\bx)-\bx)$ will be mapped to $(\sigma(\bx) \mp \bx)-(\sigma(\bx) \pm \bx) = \mp 2\bx$. We can therefore observe that if $\sigma(\bx) - \bx$ is anisotropic, then
\begin{equation}
\label{sprime1}
\sigma' = \tau_{\sigma(\bx) - \bx} \circ \sigma
\end{equation}
fixes $\bx$. If, on the other hand, $\sigma(\bx) + \bx$ is anisotropic, then
\begin{equation}
\label{sprime2}
\sigma' = \tau_{\sigma(\bx) + \bx} \circ \tau_{\sigma(\bx)} \circ \sigma
\end{equation}
fixes $\bx$. In any case, $\sigma'$ defined either by (\ref{sprime1}) or (\ref{sprime2}) is an isometry of the $(L-1)$-dimensional regular bilinear space $(\{ \bx \}^{\perp},F)$, where $\{ \bx \}^{\perp} = \{ \bz \in Z : F(\bx,\bz) = 0 \}$. Then, by the induction hypothesis,
$$\sigma' = \tau_1 \circ \dots \circ \tau_l,$$
for some reflections $\tau_1,...,\tau_l$ with $1 \leq l \leq 2L-3$ and
\begin{equation}
\label{ind}
H(\tau_i) \leq \left\{ \left( 2 N^2 |\D_K|^{\frac{1}{2d}} \right)^{\frac{(L-1)^2}{2}} H(F)^{\frac{L-1}{3}} H \left( \{ \bx \}^{\perp} \right)^{\frac{L-1}{2}} H(\sigma') \right\}^{5^{L-2}},
\end{equation}
for each $1 \leq i \leq l$, and so
\begin{equation}
\label{tau}
\sigma = \sigma'' \circ \tau_1 \circ \dots \circ \tau_l,
\end{equation}
for the same $\tau_1,...,\tau_l$ and $\sigma'' = \tau_{\sigma(\bx) - \bx}$ or $\sigma'' = \tau_{\sigma(\bx) + \bx} \circ \tau_{\sigma(\bx)}$, depending on which of $\sigma(\bx) \pm \bx$ is anisotropic, so $\sigma$ is a product of at most $2L-1$ reflections. Next we are going to produce bounds on their heights. Combining Lemma \ref{ref_height} with an argument identical to the proof of Lemma \ref{mult} and Lemma \ref{nonzero}, we obtain
\begin{equation}
\label{cd6}
H(\tau_{\sigma(\bx)}) \leq 4L N^{\frac{L+8}{2}} (N+2) |\D_K|^{\frac{L+2}{2d}} H(F) H(Z)^{\frac{L+2}{L}} H(\sigma)^2.
\end{equation}
Therefore $\tau_{\sigma(\bx)}$ satisfies (\ref{CD_me_1}). Also by Lemma \ref{ref_height},
\begin{equation}
\label{cd10}
H(\tau_{\sigma(\bx) \pm \bx}) \leq N^3(N+2) H(F) H(\sigma(\bx) \pm \bx)^2.
\end{equation}
Notice that $\sigma(\bx) \pm \bx = (A \pm I_N) \bx$. Then, once again, by an argument identical to the proof of Lemma \ref{mult}
\begin{equation}
\label{cd11}
H(\sigma(\bx) \pm \bx) \leq H(\bx) H(A \pm I_N) \leq 2 \sqrt{L} \left\{ N |\D_K|^{\frac{1}{d}} \right\}^{\frac{L+2}{4}} H(Z)^{\frac{L+2}{2L}} H(A \pm I_N),
\end{equation}
where the last inequality follows by (\ref{anis}). Combining (\ref{cd11}) with Lemma \ref{matrix}, we obtain
\begin{equation}
\label{cd12}
H(\sigma(\bx) \pm \bx) \leq 4 \sqrt{L} \left\{ N |\D_K|^{\frac{1}{d}} \right\}^{\frac{L+2}{4}} H(Z)^{\frac{L+2}{2L}} H(A).
\end{equation}
Combining (\ref{cd10}) and (\ref{cd12}), we obtain
\begin{equation}
\label{cd13}
H(\tau_{\sigma(\bx) \pm \bx}) \leq 16 L N^{\frac{L+8}{2}} (N+2) |\D_K|^{\frac{L+2}{2d}} H(F) H(Z)^{\frac{L+2}{L}} H(\sigma)^2,
\end{equation}
\smallskip
hence $\tau_{\sigma(\bx) \pm \bx}$ satisfies (\ref{CD_me_1}). By combining (\ref{sprime1}), (\ref{sprime2}), (\ref{pm}), Lemma \ref{matrix1}, (\ref{cd6}), and (\ref{cd13}), we have
\begin{equation}
\label{cd13.1}
H(\sigma') \leq 64 L^2 N^{L+8} (N+2)^2 |\D_K|^{\frac{L+2}{d}} H(F)^2 H(Z)^{\frac{2L+4}{L}} H(\sigma)^5.
\end{equation}
By Lemma \ref{intersection}, Lemma \ref{mult}, and (\ref{anis})
\begin{equation}
\label{cd13.2}
H \left( \{\bx\}^{\perp} \right) \leq H(F) H(\bx) H(Z) \leq 2 \sqrt{L} \left\{ N |\D_K|^{\frac{1}{d}} \right\}^{\frac{L+2}{4}} H(F) H(Z)^{\frac{3L+2}{2L}}.
\end{equation}
Then bound (\ref{CD_me_1}) follows upon combining (\ref{ind}) with (\ref{cd13.1}) and (\ref{cd13.2}) while keeping in mind that $2 \leq L \leq N$ and $N+2 \leq 2N$. This completes the proof.
\bigskip

{\bf Aknowledgement.} I would like to express my deep gratitude to Professor Damien Roy for pointing me in the direction of these problems, as well as for his extremely helpful suggestions that allowed to improve the bounds and simplify the arguments in this paper.

\nocite{*}
\bibliographystyle{plain}  
\bibliography{witt}        

\end{document}